\documentclass[11pt]{article}
\usepackage[T1]{fontenc}
\usepackage[utf8]{inputenc}
\usepackage{lmodern}
\usepackage[a4paper,margin=2.65cm]{geometry}
\usepackage{amsmath,amssymb,amsthm,mathtools}
\usepackage{enumitem}
\usepackage{microtype}
\usepackage{authblk}
\usepackage[colorlinks=true,linkcolor=blue,citecolor=blue,urlcolor=blue]{hyperref}

\newtheorem{thm}{Theorem}[section]
\newtheorem{prop}[thm]{Proposition}
\newtheorem{lem}[thm]{Lemma}
\newtheorem{cor}[thm]{Corollary}
\newtheorem{conj}[thm]{Conjecture}
\theoremstyle{definition}
\newtheorem{defn}[thm]{Definition}

\newcommand{\R}{\mathbb R}
\newcommand{\Z}{\mathbb Z}
\newcommand{\Q}{\mathbb Q}
\newcommand{\N}{\mathbb N}
\newcommand{\id}{\operatorname{id}}
\newcommand{\End}{\operatorname{End}}

\title{Constructive algebraic methods for bisymmetric means
}

\author[a]{Pál Burai\thanks{Email:\href{mailto:buraip@math.bme.hu}{\texttt{buraip@math.bme.hu}}.}
}
\author[b,c]{Gergely Kiss\thanks{Email:
\href{mailto:kiss.gergely@renyi.hu}{\texttt{kiss.gergely@renyi.hu}}.}}
\author[d,e]{Patricia Szokol \thanks{Email:\href{mailto:szokol.patricia@inf.unideb.hu}{\texttt{szokol.patricia@inf.unideb.hu}}.}
}
\affil[a]{Budapest University of Technology and Economics,
Egry József u. 1., 1111 Budapest, Hungary}

\affil[b]{Corvinus University of Budapest,
Fővám tér 13--15., 1093 Budapest, Hungary}
\affil[c]{HUN-REN Alfréd Rényi Institute of Mathematics,
Reáltanoda utca 13--15., 1053 Budapest, Hungary}
\affil[d]{University of Debrecen,
P.O. Box 400, H-4002 Debrecen, Hungary}
\date{}

\begin{document}
\maketitle

\begin{abstract}
We prove continuity for reflexive, bisymmetric binary operations that are strictly increasing in each variable, under certain algebraic conditions. In fact, for some algebraic weights we construct a weighted quasi-arithmetic representation of such functions without assuming continuity.
An elementary construction gives the case $r=2/3$.  More generally, a word
equality whose polynomial difference has a unique simple zero $r\in(0,1)$
forces the representation with weight $r$; every algebraic $r$ with no other
conjugate in $(0,1)$ arises this way.  The proof uses an affine embedding and a
dense-domain extension argument.
\end{abstract}


\section{Introduction}

Let $I$ be a non-degenerate interval. A binary operation $F:I^2\to I$ is \emph{reflexive} if $F(x,x)=x$, and it is
\emph{bisymmetric} if
\begin{equation}
 F(F(x,y),F(u,v))=F(F(x,u),F(y,v)).
 \label{eq:bisymmetry-intro}
\end{equation}
If $F$ is also strictly increasing in each variable, then it is a strict mean:
for $x<y$, both $F(x,y)$ and $F(y,x)$ lie strictly between $x$ and $y$.  Every
weighted quasi-arithmetic mean
\begin{equation}
 F(u,v)=\varphi^{-1}\bigl(r\varphi(u)+(1-r)\varphi(v)\bigr),
 \qquad 0<r<1,
 \label{eq:weighted-model-intro}
\end{equation}
where $\varphi$ is continuous and strictly increasing, has these three
properties.  The continuity-free problem asks for the converse.

The question originates in Acz\'el's characterization, in which continuity
is assumed in addition to the algebraic and order hypotheses
\cite{Aczel1948}.  The authors in
\cite{BuraiKissSzokol2021} removed this assumption in the symmetric binary
case.  Their subsequent dichotomy theorem
\cite{BuraiKissSzokol2023} shows that symmetry at a single pair of distinct
points already implies symmetry everywhere and hence continuity.  Thus a
reflexive, bisymmetric, strictly increasing operation is either symmetric
everywhere or satisfies $F(x,y)\ne F(y,x)$ whenever $x\ne y$.  This leaves
the nonsymmetric case as the next natural problem.

The constructive approach to removing continuity first appeared in \cite{BuraiKissSzokol2021}.
Our present approach follows the principle stated by Kiss and Shulman in \cite{KissShulman}, which consists of constructing a strictly increasing function on a dense set, establishing the required averaging identity there, and deriving continuity through a no-gap argument.

Our main question can be formalized as follows.

\begin{conj}
\label{conj:main}
Let $I\subset\R$ be a non-degenerate interval.  Every reflexive,
bisymmetric operation $F:I^2\to I$ that is strictly increasing in each
variable has the representation \eqref{eq:weighted-model-intro} for some
$r\in(0,1)$ and some continuous strictly increasing function
$\varphi:I\to\R$.
\end{conj}

This paper establishes constructive algebraic branches of
Conjecture~\ref{conj:main}.  Before stating them precisely, we briefly place
these results in their historical context.  In the regular theory of means
and functional equations, bisymmetry is one of the basic algebraic
conditions leading to quasi-arithmetic representations.  Acz\'el's early
two-variable theorem \cite{Aczel1948} was followed by work on generalized and
rectangular bisymmetry, and on characterizations with weakened or modified
hypotheses; see, among others, Acz\'el and Maksa \cite{AczelMaksa1996}, Maksa
\cite{Maksa2002}, Fodor and Marichal \cite{FodorMarichal1997}, and Marichal
\cite{Marichal2000}.  The surrounding theory of nonsymmetric and generalized
means includes the work of Bajraktarevi\'c \cite{Bajraktarevic1963}, Witold
Jarczyk and Matkowski \cite{JarczykMatkowski2000}, Matkowski
\cite{Matkowski2010}, Matkowski and P\'ales \cite{MatkowskiPales2015}, P\'ales and Pasteczka \cite{PalesPasteczka2026}
and recently,
Kiss and T\'oth \cite{KissToth2026}.

There is also a parallel algebraic history.  In universal algebra, the bisymmetry equation
\eqref{eq:bisymmetry-intro} is called the \emph{medial} or \emph{entropic law}.  The affine
structure behind medial quasigroups and related systems was developed in the
classical papers of Toyoda, Murdoch, and Bruck
\cite{Toyoda1941,Murdoch1941,Bruck1944}; Evans's abstract mean values form an
important bridge to the theory of means \cite{Evans1963}. We use this line
through the embedding theorem of Romanowska and Smith
\cite{RomanowskaSmith2001} stated in Section~2; see also Stronkowski
\cite{StronkowskiModules}.

The dense-domain extension step is stated as
Theorem~\ref{thm:dense-domain} below, in the binary form of
\cite[Proposition~7.1]{KissShulman}.  This step requires neither symmetry
nor bisymmetry.  In the present paper bisymmetry and reflexivity are used
earlier, to construct the ordered weighted structure to which the theorem
can be applied.

In parallel with this work, Kiss and Okamura \cite{KissOkamura2026} obtained an existential proof
of the characterization of weighted quasi-arithmetic means.  We present
here the results of the constructive approach, which explicitly builds the
representation directly from algebraic word information.

Returning to the constructive problem, we introduce only the notation needed
to state the results transparently. Write $u*v:=F(u,v)$. One of the simplest
nonsymmetric branches is governed by one depth-two coincidence and permits a
completely elementary construction.

\begin{thm}
\label{thm:intro-23}
Let $F:I^2\to I$ be reflexive, bisymmetric, and strictly increasing in each
variable.  If $x<y$ and
\begin{equation}
 (y*x)*x=(x*y)*(y*x),
 \label{eq:model-identity-intro}
\end{equation}
then there is a continuous strictly increasing bijection
$\varphi:[x,y]\to[0,1]$ such that
\begin{equation}
 \varphi(F(u,v))=\frac23\varphi(u)+\frac13\varphi(v)
 \qquad(u,v\in[x,y]).
 \label{eq:23-result-intro}
\end{equation}
Moreover, there is a continuous strictly increasing function
$\psi:I\to\R$ such that
\[
 F(u,v)=\psi^{-1}\left(
 \frac23\psi(u)+\frac13\psi(v)
 \right)
 \qquad(u,v\in I).
\]
\end{thm}

For the affine operation $u*v=ru+(1-r)v$, the identity \eqref{eq:model-identity-intro}, evaluated  at $x=0$, $y=1$, becomes $r^2=2r(1-r)$.
Thus, the identity selects the weight $r=2/3$ among $0<r<1$.
Theorem~\ref{thm:intro-23} proves that a coincidence at one pair $x<y$
already forces the corresponding weighted representation on the whole
interval $I$, without assuming continuity.

To formulate a general algebraic branch we need the following concepts.
Let $I$ be a non-degenerate real interval and $*$ a binary operation on $I$, that is,
\[
 *\colon I\times I\to I.
\]
Let $x,y\in I$ be arbitrarily fixed. We construct binary words in two
formal letters $X,Y$, representing $x,y$, respectively, in the following
recursive way.
\begin{itemize}
\item The zeroth-level binary words are
\[
 X,\qquad Y.
\]
\item The first-level binary words are
\[
 X*X,\qquad X*Y,\qquad Y*X,\qquad Y*Y.
\]
\item If $U,V$ are binary words of levels $k,\ell$, respectively, then
\[
 U*V
\]
is a binary word of level $\max\{k,\ell\}+1$. Every binary word of
positive level is constructed in this way.
\end{itemize}

Let us denote by $\mathcal{W}_*(X,Y)$ the set of all binary words
constructed from formal letters $X$ and $Y$. If there is no ambiguity, we just write
$\mathcal{W}$ instead of $\mathcal{W}_*(X,Y)$.
Evaluation at $x,y$ is the map
\[
 \mathcal{W}\longrightarrow I,\qquad W\longmapsto W(x,y),
\]
obtained by substituting $x,y$ for $X,Y$ and applying $*$.
This map need not be injective. The polynomials below are assigned to
the formal words, before evaluation.

We assign a polynomial to every word in $\mathcal{W}$ in the following
recursive way. Let
\[
 P_X(t):=0,\qquad P_Y(t):=1,
\]
and, for every $U,V\in\mathcal{W}$, let
\[
 P_{U*V}(t):=tP_U(t)+(1-t)P_V(t).
\]

In the weighted quasi-arithmetic case \eqref{eq:weighted-model-intro}, the same recursion gives
\[
 \varphi(W(x,y))
 =(1-P_W(r))\varphi(x)+P_W(r)\varphi(y).
\]
The following theorem identifies the word coincidences treated by our
method and shows that they realize every algebraic weight whose conjugates
do not create a second admissible weight.

\begin{thm}
\label{thm:intro-main}
Let $F:I^2\to I$ be reflexive, bisymmetric, and strictly increasing in each
variable.  Let $U,V$ be two different binary words and suppose that the nonzero polynomial
\[
 Q=P_{U}-P_{V}
\]
has exactly one zero $r$ in $(0,1)$ and that this zero is simple.  Let $a<b$ be points of $I$ such that
$U(a,b)=V(a,b)$.
Then the word identity
\[U(x,y)=V(x,y)\] holds for every $x,y\in I$ and
there is a continuous strictly increasing function $\varphi:I\to\R$ such that
\[
 F(u,v)=\varphi^{-1}\bigl(r\varphi(u)+(1-r)\varphi(v)\bigr)
 \qquad(u,v\in I).
\]

Conversely, let $r\in(0,1)$ be algebraic and suppose that its minimal
polynomial over $\Q$ has no other zero in $(0,1)$.  Then there are binary words
$U_r,V_r$ for which $P_{U_r}-P_{V_r}$ is a polynomial multiple of that
minimal polynomial and has $r$ as its unique simple zero in $(0,1)$.
\end{thm}

The proof constructs a strictly increasing function
$f_r:S_r\to[a,b]$ by
\[
 S_r:=\{P_W(r):W\in\mathcal{W}\},\qquad
 f_r(P_W(r)):=W(a,b).
\]
The polynomial arguments establish that this assignment is well-defined
and strictly increasing, and that $S_r$ is countable and dense in $[0,1]$.
Its exact weighted averaging identity permits the application of
Theorem~\ref{thm:dense-domain}
\cite[Proposition~7.1]{KissShulman}, which proves density of the image
$f_r(S_r)$ in $[a,b]$ and gives the continuous local representation.
This local conclusion also shows that the equivalence classes associated
with the word equality are intervals. The globalization principle in
Proposition~\ref{prop:global-propagation} and
Corollary~\ref{cor:polynomial-globalization} then gives the identity and
the representation on the whole interval.

The unique-root hypothesis also marks the scope of the result. One rational
word equality cannot distinguish several interior conjugates, and a
transcendental weight cannot satisfy a nonzero polynomial equality.
Section~6 discusses the additional information that may be needed beyond
this case.

The structure of the paper is as follows.
Section~2 develops the affine representation and word evaluation,
then treats polynomial relations and globalization, and ends with the
dense-domain completion theorem. Section~3 treats the weight $2/3$ (Theorem \ref{thm:intro-23}).
Section~4 establishes the polynomial tools, with explicit finite proofs of
the positive-coefficient statements. Section~5 proves the unique-root
theorem (Theorem~\ref{thm:intro-main}), classifies depth-two word coincidences, and gives a family of higher-depth examples. Section~6 discusses further directions.
The two order-theoretic interval arguments are proved in Appendices~A and~B.

\section{Preliminaries}

Throughout the paper $I\subset\R$ is a non-degenerate interval and
$F:I^2\to I$ is reflexive, bisymmetric and strictly increasing in each
variable.  With $x*y:=F(x,y)$, these assumptions are
\begin{align}
 x*x&=x,                                                   \label{eq:idempotent}\\
 (x*y)*(u*v)&=(x*u)*(y*v),                                \label{eq:medial}\\
 x<x'&\Longrightarrow x*y<x'*y,\qquad
 y<y'\Longrightarrow x*y<x*y'.                           \label{eq:strict}
\end{align}
A set equipped with a reflexive and bisymmetric binary operation\footnote{In universal algebra, reflexivity is also called idempotency, and bisymmetry is called the medial or entropic law.} is called a
\emph{binary mode}. It is \emph{cancellative} if each map $u\mapsto u*z$
and $u\mapsto z*u$ is injective. Thus \eqref{eq:idempotent} and
\eqref{eq:medial} make $(I,*)$ a binary mode, and strict increase implies
cancellativity.

\begin{lem}
\label{lem:strict-mean}
If $u<v$, then
\[
 u<u*v<v,
 \qquad
 u<v*u<v.
\]
Consequently, every compact subinterval of $I$ is closed under $*$.
\end{lem}

\begin{proof}
By reflexivity and strict increase,
\[
 u=u*u<u*v<v*v=v,
\]
and the second inequality is analogous.
\end{proof}

\subsection{Affine representation and word polynomials}

We use the following binary form of the affine embedding theorem of
Romanowska and Smith
\cite[Theorem~3.5 and Corollary~3.6]{RomanowskaSmith2001}.

\begin{thm}
\label{thm:affine-input}
Let $*$ be a reflexive and bisymmetric binary operation on a nonempty
set $A$. Suppose that $*$ is cancellative, i.e., for all $u,v,z\in A$,
\[
 u*z=v*z\ \Longrightarrow\ u=v,
 \qquad
 z*u=z*v\ \Longrightarrow\ u=v.
\]
Then there exist an abelian group $M$, additive bijections
$\alpha,\beta:M\to M$, and an injective map $\iota:A\to M$ such that
\[
 \alpha\circ\beta=\beta\circ\alpha,
 \qquad
 \alpha+\beta=\id,
\]
and
\begin{equation}
 \iota(u*v)=\alpha(\iota(u))+\beta(\iota(v))
 \qquad(u,v\in A).
 \label{eq:affine-input}
\end{equation}
\end{thm}

Here $\alpha+\beta=\id$ means that
$\alpha(w)+\beta(w)=w$ for every $w\in M$.
Both maps are injective, so $\alpha(p)=\alpha(q)$ implies $p=q$,
and similarly for $\beta$. These cancellation properties will be
used below. Strict increase ensures that the theorem applies to
our operation on $I$. We fix such an embedding of $(I,*)$ for the rest
of the paper and represent each subalgebra by restriction.

We first explain how the embedding is used. For a chosen base point
$a\in I$, let
\[
 j(u):=\iota(u)-\iota(a).
\]
Thus $j(a)=0$; for a second point $b$ we write
$d:=j(b)=\iota(b)-\iota(a)$. Translating the image does not alter the
formula for the operation. Indeed,
\begin{align*}
 j(u*v)
 &=\alpha(\iota(u))+\beta(\iota(v))-\iota(a)\\
 &=\alpha(j(u))+\beta(j(v))
   +(\alpha+\beta-\id)(\iota(a))\\
 &=\alpha(j(u))+\beta(j(v)).
\end{align*}
Consequently we shall use
\begin{equation}
 j(u*v)=\alpha(j(u))+\beta(j(v)),\qquad
 \beta=\id-\alpha,\qquad \alpha\beta=\beta\alpha.
 \label{eq:affine-notation}
\end{equation}

The symbols $\alpha,\beta$ denote additive maps on $M$.
Products of maps and
powers denote composition, for example
$(\alpha\beta)(d)=\alpha(\beta(d))$ and $\alpha^2(d)=\alpha(\alpha(d))$.
For $P(t)=\sum_k p_kt^k\in\Z[t]$, we define
\[
 P(\alpha)(d):=\sum_k p_k\alpha^k(d),\qquad \alpha^0=\id.
\]
The integer multiples and sums are taken in the abelian group.
Polynomial calculations are legitimate because the commuting maps
generate a commutative subring of $\End(M)$.

We keep the embedding explicit: a difference of images is written
$j(v)-j(u)$, while an inequality $u<v$ refers to the original real points.
No order or topology is assumed on $M$. In particular, $j(a)=0$ and
$j(b)=d$ are statements about their images; they do not replace the real
points $a,b$ by $0,d$. A continuous real-valued function representing $F$
will be constructed from the ordered word values later.

\begin{defn}
For the formal binary words in $X,Y$ introduced above, evaluation at
$u,v\in I$ is defined by $X(u,v)=u$, $Y(u,v)=v$, and
$(U*V)(u,v)=U(u,v)*V(u,v)$. The word polynomials satisfy
\begin{equation}
 P_X(t)=0,
 \qquad P_Y(t)=1,
 \qquad
 P_{U*V}(t)=tP_U(t)+(1-t)P_V(t).
 \label{eq:word-polynomial}
\end{equation}
Thus $P_U\in\Z[t]$ for every word $U$.
\end{defn}

\begin{lem}
\label{lem:word-evaluation}
Let $a,b\in I$, use the fixed affine embedding, and let
\[
 j(u)=\iota(u)-\iota(a),\qquad d=j(b).
\]
Then, for every binary word $U$,
\begin{equation}
 j(U(a,b))=P_U(\alpha)(d).
 \label{eq:word-evaluation}
\end{equation}
Equivalently, $\iota(U(a,b))=\iota(a)+P_U(\alpha)(d)$.
\end{lem}

\begin{proof}
For the two generators, the formula says $j(a)=0$ and $j(b)=d$.
If it holds for $U,V$, then \eqref{eq:affine-notation} gives
\[
 j((U*V)(a,b))
 =\alpha(P_U(\alpha)(d))+\beta(P_V(\alpha)(d))
 =P_{U*V}(\alpha)(d),
\]
which proves the induction step.
\end{proof}

The same formula, with $u$ as base point, gives for arbitrary $u,v\in I$
and every binary word $W$,
\begin{equation}
 j(W(u,v))=j(u)+P_W(\alpha)\bigl(j(v)-j(u)\bigr).
 \label{eq:word-evaluation-any-pair}
\end{equation}
This identity is independent of the base point used to define $j$.

For example, Lemma~\ref{lem:word-evaluation} gives
\[
 j(a*b)=\beta(d),\qquad j(b*a)=\alpha(d),\qquad
 j((b*a)*a)=\alpha^2(d).
\]
These formulas turn word equalities into relations between endomorphisms
applied to $d$. Recovering the order of word values requires the strict
increase of the original operation.

\subsection{Word equalities and globalization}

We first collect the algebraic consequences of a word coincidence.
For an equivalence relation $\sim$ on an interval closed under $*$, we
use the following condition for all triples of points in that interval:
\begin{equation}
 x*z\sim y*z\ \Longleftrightarrow\ x\sim y,
 \qquad
 z*x\sim z*y\ \Longleftrightarrow\ x\sim y.
 \label{eq:quotient-cancellation}
\end{equation}
These equivalences express both preservation of classes and cancellation
when one argument is fixed.

\begin{lem}
\label{lem:word-coincidence-subalgebra}
Use the fixed affine embedding of $(I,*)$.
\begin{enumerate}[label=\textup{(\roman*)},leftmargin=*]
\item For $P\in\Z[t]$, define
\[
 x\sim_P y\quad\Longleftrightarrow\quad
 P(\alpha)(j(y)-j(x))=0.
\]
Then $\sim_P$ is an equivalence relation satisfying
\eqref{eq:quotient-cancellation}.
\item Let $U,V$ be binary words and put $Q=P_U-P_V$. Then
\[
 x\sim_Q y\quad\Longleftrightarrow\quad U(x,y)=V(x,y).
\]
If $U(a,b)=V(a,b)$ for some $a<b$, the common equivalence class of
$a,b$ contains the subalgebra
$A:=\{W(a,b):W\in\mathcal W\}$. Consequently,
\[
 U(u,v)=V(u,v)\qquad(u,v\in A).
\]
\end{enumerate}
\end{lem}

\begin{proof}
\textup{(i)} Put $T=P(\alpha)$. Additivity gives
\[
 x\sim_P y\quad\Longleftrightarrow\quad T(j(x))=T(j(y)),
\]
so $\sim_P$ is an equivalence relation.
Since $T$ commutes with $\alpha,\beta$, the affine formula gives
\begin{align*}
 T(j(y*z)-j(x*z))&=\alpha\bigl(T(j(y)-j(x))\bigr),\\
 T(j(z*y)-j(z*x))&=\beta\bigl(T(j(y)-j(x))\bigr).
\end{align*}
Injectivity of $\alpha,\beta$ now proves both equivalences in
\eqref{eq:quotient-cancellation}.

\textup{(ii)} The equivalence between the word equality and $\sim_Q$
follows from \eqref{eq:word-evaluation-any-pair} and injectivity of $j$.
Let $C$ be the class containing $a,b$. If $u,v\in C$, part~\textup{(i)}
and reflexivity give
\[
 u*v\sim_Q a*v\sim_Q a*a=a.
\]
Thus $C$ is closed under $*$ and contains every word value generated by
$a,b$. Hence all $u,v\in A$ satisfy $u\sim_Q v$, which proves the
required word identity on $A$.
\end{proof}

The following principle adapts the equivalence-class strategy of
Burai, Kiss, and Szokol
\cite[Section~3.2, proof of Theorem~3]{BuraiKissSzokol2023}.
There the relation $x\sim y$ defined by $x*y=y*x$ propagates symmetry
from one pair of distinct points to the whole interval. The interval
property of its classes is first established by the local representation
in \cite[Lemma~6]{BuraiKissSzokol2023}. Here we isolate the subsequent
interval argument for application to polynomial relations.
In the following statements, equivalence classes are required to be
intervals, possibly singletons. Equivalently, whenever two points belong
to the same class, so does every point between them.

\begin{prop}
\label{prop:global-propagation}
Let $K$ be a non-degenerate compact real interval and let $*:K^2\to K$
be reflexive and strictly increasing in each variable. Let $\sim$ be an
equivalence relation whose classes are intervals. Suppose that
\eqref{eq:quotient-cancellation} holds for all $x,y,z\in K$.
If one class contains two distinct points, then $K$ consists of a single
class.
\end{prop}

\begin{proof}
The proof is given in Appendix~B. The role of
\eqref{eq:quotient-cancellation} is that translations preserve classes and
keep distinct classes separate. A class containing an interval would then
produce a disjoint interval for each distinct translated class. The proof
combines this observation with compactness of $K$.
\end{proof}

Here is the form in which we shall apply the interval principle to the
polynomial relations from Lemma~\ref{lem:word-coincidence-subalgebra}.

\begin{cor}
\label{cor:polynomial-globalization}
Let $P\in\Z[t]$, and let $\sim_P$ be the relation from
Lemma~\ref{lem:word-coincidence-subalgebra}(i).
Suppose its classes are intervals. If $a\sim_P b$ for some $a<b$, then
\[
 P(\alpha)(j(y)-j(x))=0\qquad(x,y\in I).
\]
\end{cor}

\begin{proof}
By Lemma~\ref{lem:word-coincidence-subalgebra}(i), the relation satisfies
\eqref{eq:quotient-cancellation}. Fix $x,y\in I$, and choose a compact
subinterval $K\subset I$ containing $a,b,x,y$.
Lemma~\ref{lem:strict-mean} ensures that $K$ is closed under $*$.
The restricted relation still has interval classes and satisfies
\eqref{eq:quotient-cancellation}. Since $a<b$ belong to the same class,
Proposition~\ref{prop:global-propagation} shows that $K$ consists of one
class. In particular, $x\sim_P y$.
\end{proof}

What remains in applications is to verify the interval property:
whenever $u<v$ and $u\sim_P v$, every point of $[u,v]$ must lie in their
class. For relations arising from word equalities,
Lemma~\ref{lem:word-coincidence-subalgebra}(ii) ensures that their class
contains the subalgebra generated by $u,v$.
The stronger interval property will follow from the local
representations in Corollary~\ref{cor:global-23} and in Step~4 of the
proof of Theorem~\ref{thm:unique-root}.

\subsection{Dense-domain completion}

The local constructions in Sections~3 and~5 first produce a strictly
increasing function on a dense parameter set, together with an exact
weighted identity. The following theorem extends that function to the
whole compact interval and gives the local representation used to verify
the interval hypothesis above.

\begin{thm}
\label{thm:dense-domain}
Let $J=[a,b]$ be a non-degenerate compact interval, let $0<r<1$, and let
$H:J^2\to J$ be strictly increasing in each variable.  Suppose that
$S\subset[0,1]$ is dense, contains $0,1$, and is closed under
\[
 (s,t)\longmapsto rs+(1-r)t.
\]
If $f:S\to J$ is strictly increasing, $f(0)=a$, $f(1)=b$, and
\begin{equation}
 H(f(s),f(t))=f\bigl(rs+(1-r)t\bigr)
 \qquad(s,t\in S),
 \label{eq:dense-domain-identity}
\end{equation}
then $f(S)$ is dense in $J$.  Consequently $f$ extends uniquely to a
continuous strictly increasing bijection $\widetilde f:[0,1]\to J$, and
\begin{equation}
 H(u,v)=\widetilde f\bigl(
 r\widetilde f^{-1}(u)+(1-r)\widetilde f^{-1}(v)
 \bigr)
 \qquad(u,v\in J).
 \label{eq:dense-domain-conclusion}
\end{equation}
\end{thm}

\begin{proof}
The no-gap idea already appears in the symmetric binary argument of
\cite{BuraiKissSzokol2021}; the statement above is the binary specialization
of the general dense-domain theorem in
\cite[Proposition~7.1]{KissShulman}.  It uses strict monotonicity, but neither
symmetry nor bisymmetry.  For completeness, a short proof specialized to the
present binary setting is included in Appendix~A.
\end{proof}

\section{The constructive \texorpdfstring{$2/3$}{2/3} branch}

Fix $x<y$ and assume
\begin{equation}
 (y*x)*x=(x*y)*(y*x).
 \label{eq:23-identity}
\end{equation}

\subsection{The first subdivision}

\begin{lem}
\label{lem:first-order}
If $u<v$ and
\begin{equation}
 (v*u)*u=(u*v)*(v*u),
 \label{eq:local-23-identity}
\end{equation}
then
\[
 u<u*v<v*u<v.
\]
\end{lem}

\begin{proof}
Lemma~\ref{lem:strict-mean} places both $u*v$ and $v*u$ strictly between
$u$ and $v$. Put $p=u*v$ and $q=v*u$. The identity says $q*u=p*q$.
If $q=p$, then
\[
 q*u=p*u<p=p*p=p*q,
\]
a contradiction. If $q<p$, then
\[
 q*u<q<p*q,
\]
again a contradiction. Hence $p<q$.
\end{proof}

The subdivision construction will repeatedly use
\eqref{eq:23-identity} at pairs obtained from $x,y$.
Let $A:=\{W(x,y):W\in\mathcal W\}$. By
Lemma~\ref{lem:word-coincidence-subalgebra}(ii), the initial coincidence gives
\[
 (v*u)*u=(u*v)*(v*u)\qquad(u,v\in A).
\]
We now express the coincidence as a linear relation in the affine group;
this relation will control the adjacent differences in the subdivisions.

Restrict the fixed affine embedding to $A$, and take $x$ as the base point.
Thus $j(x)=0$ and $d=j(y)$.
We keep $j$ explicit throughout the calculation, so
\[
 j(u*v)=\alpha(j(u))+\beta(j(v)),
 \qquad \alpha+\beta=\id.
\]

We have
\[
 j(y*x)=\alpha(d),
 \qquad
 j((y*x)*x)=\alpha(\alpha(d))+\beta(0)=\alpha^2(d),
\]
and $j(x*y)=\beta(d)$. Hence, using the commutativity of
$\alpha$ and $\beta$,
\begin{align*}
 j((x*y)*(y*x))
 &=\alpha(\beta(d))+\beta(\alpha(d))\\
 &=2\alpha(\beta(d)).
\end{align*}
Consequently \eqref{eq:23-identity} gives
\[
 \alpha^2(d)=2\alpha(\beta(d)),
 \qquad\text{that is,}\qquad
 \alpha\bigl(\alpha(d)-2\beta(d)\bigr)=0.
\]
Since $\alpha$ is injective, we obtain
\[
 \alpha(d)-2\beta(d)=0.
\]
As $\beta(d)=d-\alpha(d)$, this is equivalent to
\begin{equation}
 (3\alpha-2\id)(d)=0.
 \label{eq:linear-relation}
\end{equation}

The following proposition gives the ordered subdivision of any pair
$u<v$ in $A$ and verifies \eqref{eq:linear-relation} for the three
new adjacent differences.

\begin{prop}
\label{prop:propagation}
Let $u<v$ belong to $A$. Then
\[
 u<u*v<v*u<v,
\]
and $3\alpha-2\id$ annihilates each of the three adjacent image differences
\[
 j(u*v)-j(u),\qquad
 j(v*u)-j(u*v),\qquad
 j(v)-j(v*u).
\]
\end{prop}

\begin{proof}
Lemma~\ref{lem:word-coincidence-subalgebra}(ii) gives
\eqref{eq:local-23-identity} for $u,v$, so the ordering follows from
Lemma~\ref{lem:first-order}. Put $e=j(v)-j(u)$.
The word-polynomial difference in \eqref{eq:local-23-identity} is
$t(3t-2)$. Thus \eqref{eq:word-evaluation-any-pair} gives
$\alpha\bigl((3\alpha-2\id)(e)\bigr)=0$. Injectivity of $\alpha$ yields
$(3\alpha-2\id)(e)=0$. Moreover,
\begin{align*}
 j(u*v)-j(u)&=\beta(e),\\
 j(v*u)-j(u*v)&=(\alpha-\beta)(e),\\
 j(v)-j(v*u)&=\beta(e).
\end{align*}
Since $3\alpha-2\id$ commutes with $\alpha$ and $\beta$, it annihilates
all three displayed differences.
\end{proof}

\subsection{The triadic coding}

Proposition~\ref{prop:propagation} supplies the local step of the recursive
construction: each adjacent pair is replaced by four correctly ordered
points, and the same relation is inherited by the three new adjacent pairs.
All new points remain in $A$, so the proposition can be applied at every
stage. This produces the canonical ordered representatives of the triadic numbers.

For $n\ge0$ put
\[
 D_{3,n}:=\left\{\frac{k}{3^n}:0\le k\le3^n\right\},
 \qquad
 D_3:=\bigcup_{n\ge0}D_{3,n}.
\]
We recursively define formal words $W_{n,k}$ in $X,Y$.  Set
\[
 W_{0,0}=X,
 \qquad W_{0,1}=Y.
\]
If $0\le k<3^n$, write $L=W_{n,k}$ and $U=W_{n,k+1}$ and put
\begin{align}
 W_{n+1,3k}&=L,&
 W_{n+1,3k+1}&=L*U,\notag\\
 W_{n+1,3k+2}&=U*L,&
 W_{n+1,3k+3}&=U.
 \label{eq:triadic-recursion}
\end{align}
The two definitions at a common endpoint of neighboring cells agree,
since both retain the same old word. Define
\begin{equation}
 f_n\left(\frac{k}{3^n}\right):=W_{n,k}(x,y).
 \label{eq:fn}
\end{equation}

\begin{prop}
\label{prop:fn-increasing}
For every $n$, the map $f_n:D_{3,n}\to[x,y]$ is strictly increasing,
$f_{n+1}$ extends $f_n$, and every adjacent difference at level $n$ is
annihilated by $3\alpha-2\id$, where the adjacent difference means
$j(f_n((k+1)/3^n))-j(f_n(k/3^n))$. Hence
\[
 f:=\bigcup_{n\ge0}f_n:D_3\longrightarrow[x,y]
\]
is well-defined and strictly increasing.
\end{prop}

\begin{proof}
At level zero, $f_0(0)=x<y=f_0(1)$, and the operator relation is
\eqref{eq:linear-relation}. If $L<U$ are adjacent
values at level $n$, Proposition~\ref{prop:propagation} gives
\[
 L<L*U<U*L<U
\]
and propagates the operator relation to the three new adjacent differences.
The recursion retains every old word, so $f_{n+1}$ extends $f_n$.
The assertion follows by induction.
\end{proof}

\begin{lem}
\label{lem:canonical-coordinate}
For $0\le k\le3^n$,
\begin{equation}
 P_{W_{n,k}}(2/3)=\frac{k}{3^n}.
 \label{eq:canonical-coordinate}
\end{equation}
Moreover,
\begin{equation}
 \{P_U(2/3):U\text{ a binary word}\}=D_3.
 \label{eq:coordinate-set}
\end{equation}
\end{lem}

\begin{proof}
The first assertion follows by induction from
\eqref{eq:triadic-recursion}.  Indeed, the two new points between
$k/3^n$ and $(k+1)/3^n$ are
\[
 \frac23\frac{k}{3^n}+\frac13\frac{k+1}{3^n}
 =\frac{3k+1}{3^{n+1}}
\]
and
\[
 \frac23\frac{k+1}{3^n}+\frac13\frac{k}{3^n}
 =\frac{3k+2}{3^{n+1}}.
\]
Thus every triadic rational occurs.  Conversely, $D_3$ contains $0,1$ and is
closed under $(s,t)\mapsto(2s+t)/3$, so induction on words gives the reverse
inclusion in \eqref{eq:coordinate-set}.
\end{proof}

\subsection{Consistency of the word coding}

The preceding recursion defines one canonical word for each triadic number.
For the weighted identity, however, we must be able to use arbitrary words
$U,V$ representing arbitrary $s,t\in D_3$.  We therefore have to prove that
two words with the same polynomial value at $2/3$ always have the same value at
$(x,y)$.  The following elementary form of Gauss's lemma is the algebraic core
of this consistency statement.

\begin{lem}
\label{lem:rational-root-divisibility}
If $Q\in\Z[t]$ and $Q(2/3)=0$, then there is $R\in\Z[t]$ such that
\[
 Q(t)=(3t-2)R(t).
\]
\end{lem}

\begin{proof}
The polynomial $3t-2$ is primitive and irreducible in $\Z[t]$.  It divides
$Q$ in $\Q[t]$, and Gauss's lemma therefore gives divisibility already in
$\Z[t]$.
\end{proof}

\begin{prop}
\label{prop:word-consistency}
For arbitrary binary words $U,V$,
\[
 P_U(2/3)=P_V(2/3)
 \quad\Longrightarrow\quad
 U(x,y)=V(x,y).
\]
Consequently
\begin{equation}
 f(P_U(2/3))=U(x,y)
 \label{eq:consistent-assignment}
\end{equation}
for every binary word $U$.
\end{prop}

\begin{proof}
Put $Q=P_U-P_V$.  By Lemma~\ref{lem:rational-root-divisibility},
$Q=(3t-2)R$ for some $R\in\Z[t]$.  Hence
\[
 Q(\alpha)(d)
 =R(\alpha)\bigl((3\alpha-2\id)(d)\bigr)=0
\]
by \eqref{eq:linear-relation}.  Lemma~\ref{lem:word-evaluation} now gives
$U(x,y)=V(x,y)$.

For \eqref{eq:consistent-assignment}, choose $n,k$ with
$P_U(2/3)=k/3^n$.  Lemma~\ref{lem:canonical-coordinate} gives a canonical word
$W_{n,k}$ with the same polynomial value at $2/3$, and therefore
\[
 U(x,y)=W_{n,k}(x,y)=f(k/3^n).
\]
\end{proof}

\subsection{The exact weighted identity on the dense set}

Proposition~\ref{prop:fn-increasing} gives an ordered canonical coding, while
Proposition~\ref{prop:word-consistency} makes it independent of the chosen
word.  The latter is precisely what permits the operation to be evaluated at
an arbitrary pair of points of $D_3$.

\begin{prop}
\label{prop:dense-identity}
For all $s,t\in D_3$,
\begin{equation}
 F(f(s),f(t))=f\left(\frac{2s+t}{3}\right).
 \label{eq:23-dense-identity}
\end{equation}
\end{prop}

\begin{proof}
Choose words $U,V$ with $s=P_U(2/3)$ and $t=P_V(2/3)$.  By
Proposition~\ref{prop:word-consistency}, $f(s)=U(x,y)$ and
$f(t)=V(x,y)$.  Since
\[
 P_{U*V}(2/3)=\frac23P_U(2/3)+\frac13P_V(2/3),
\]
another application of Proposition~\ref{prop:word-consistency} gives
\[
 f\left(\frac{2s+t}{3}\right)
 =(U*V)(x,y)=F(f(s),f(t)).
\]
\end{proof}

\begin{thm}
\label{thm:local-23}
Let $F:I^2\to I$ be reflexive, bisymmetric and strictly increasing in each
variable. Suppose that $x<y$ and
\[
 (y*x)*x=(x*y)*(y*x).
\]
Then there is a continuous strictly increasing bijection
$\varphi:[x,y]\to[0,1]$ such that
\begin{equation}
 \varphi(F(u,v))=\frac23\varphi(u)+\frac13\varphi(v)
 \qquad(u,v\in[x,y]).
 \label{eq:local-23-representation}
\end{equation}
In particular, $F$ is continuous on $[x,y]^2$.
\end{thm}

\begin{proof}
The set $D_3$ is dense in $[0,1]$, contains $0,1$, and is closed under
$(s,t)\mapsto(2s+t)/3$.  Proposition~\ref{prop:fn-increasing} says that
$f:D_3\to[x,y]$ is strictly increasing, while
Proposition~\ref{prop:dense-identity} gives the exact weighted identity on
$D_3$.  The endpoint values are $f(0)=x$ and $f(1)=y$.  Thus all hypotheses
of Theorem~\ref{thm:dense-domain}
\cite[Proposition~7.1]{KissShulman} hold with $H=F|_{[x,y]^2}$ and $r=2/3$.  If
$\widetilde f:[0,1]\to[x,y]$ denotes the resulting extension, put
$\varphi=\widetilde f^{-1}$.
\end{proof}

\begin{cor}
\label{cor:global-23}
Suppose that, for some $a<b$ in $I$,
\begin{equation}
 (b*a)*a=(a*b)*(b*a).
 \label{eq:global-23-identity}
\end{equation}
Then the same identity holds for every pair of points of $I$.
Moreover, $F$ is continuous on $I^2$ and there is a
continuous strictly increasing function $\varphi:I\to\R$, unique up to a
positive affine change, such that
\[
 F(u,v)=\varphi^{-1}\left(
 \frac23\varphi(u)+\frac13\varphi(v)
 \right)
 \qquad(u,v\in I).
\]
\end{cor}

\begin{proof}
Define a relation on $I$ by
\[
 u\sim v\quad\Longleftrightarrow\quad
 (v*u)*u=(u*v)*(v*u).
\]
By Lemma~\ref{lem:word-coincidence-subalgebra}(i)--(ii), this is the
relation $\sim_P$ for $P(t)=t(3t-2)$. It is therefore an equivalence
relation satisfying \eqref{eq:quotient-cancellation}.

Let $u<v$ with $u\sim v$. Theorem~\ref{thm:local-23} gives a continuous
strictly increasing bijection $\theta:[u,v]\to[0,1]$ representing $F$
with weight $2/3$. For every $s,t\in[u,v]$, this representation gives
\[
 \theta((t*s)*s)
 =\frac59\theta(s)+\frac49\theta(t)
 =\theta((s*t)*(t*s)).
\]
Since $\theta$ is injective, $s\sim t$. Hence all points of $[u,v]$
lie in one equivalence class, so each class is an interval.
The assumed coincidence gives $a\sim b$ with $a<b$.
Corollary~\ref{cor:polynomial-globalization} therefore gives $u\sim v$
for all $u,v\in I$, which proves the global word identity.

Each point of $I^2$ has a relative neighborhood contained in $K^2$ for
some non-degenerate compact subinterval $K\subset I$. By
Theorem~\ref{thm:local-23}, $F$ is continuous on $K^2$, and hence on $I^2$.
The classical continuous bisymmetric representation theorem gives a
weighted quasi-arithmetic form with some weight $r\in(0,1)$
\cite[Section~2]{Aczel1948}. Substitution into
\eqref{eq:global-23-identity} gives $r^2=2r(1-r)$, and therefore $r=2/3$.

For uniqueness, let $\varphi_1,\varphi_2$ be two increasing continuous
generators and put $h=\varphi_2\circ\varphi_1^{-1}$. Then
$h((2s+t)/3)=(2h(s)+h(t))/3$. Applying the canonical words from
Lemma~\ref{lem:canonical-coordinate} to any two points $s<t$ shows
\[
 h((1-q)s+qt)=(1-q)h(s)+qh(t)\qquad(q\in D_3).
\]
Continuity extends this equality to every $q\in[0,1]$. Thus $h$ is affine
on every compact subinterval of its domain. Fixing two distinct points
determines the same affine formula on every compact subinterval containing
them, and such subintervals cover the domain. Therefore $h(s)=cs+c_0$
with $c>0$.
\end{proof}

\begin{proof}[Proof of Theorem~\ref{thm:intro-23}]
The local assertion is Theorem~\ref{thm:local-23}, and the global assertion
is Corollary~\ref{cor:global-23}, applied to the same pair $x<y$.
\end{proof}

\section{Polynomial tools}

The general algebraic case needs three facts before the main proof begins:
integer polynomials can be realized by differences of word polynomials;
strictly positive polynomials act as positive combinations of translated
word differences; and the Euclidean algorithm can be used to choose
a polynomial multiplier positive on $[0,1]$. All three statements are
proved below by finite calculations.

\subsection{Complete words and polynomial differences}

Complete a shorter word to a full binary expression of depth $n$ by
replacing leaves $Z$ with $Z*Z$. Reflexivity preserves its value, and
$tP_Z+(1-t)P_Z=P_Z$ preserves its polynomial. A path to a leaf uses the
first argument $k$ times and the second argument $n-k$ times for some $k$.
There are $\binom nk$ such paths. A leaf labelled $Y$ contributes
$t^k(1-t)^{n-k}$, while a leaf labelled $X$ contributes zero. Therefore
the polynomials of complete words of depth $n$ are exactly
\begin{equation}
 P(t)=\sum_{k=0}^n c_k t^k(1-t)^{n-k},\qquad
 c_k\in\Z,\quad 0\le c_k\le\binom nk.
 \label{eq:Bernstein-word}
\end{equation}
Indeed, for each $k$ we may independently label any $c_k$ of the
$\binom nk$ corresponding leaves by $Y$.

\begin{lem}
\label{lem:polynomial-realization}
For every $Q\in\Z[t]$ there are an integer $N\ge0$ and binary words
$U,V$ of a common depth such that
\begin{equation}
 P_U(t)-P_V(t)=t^N(1-t)^NQ(t).
 \label{eq:polynomial-realization}
\end{equation}
\end{lem}

\begin{proof}
For $Q=0$ take $U=V=X$ and $N=0$. Otherwise write
$Q(t)=\sum_{j=0}^q a_jt^j$ and put $n=2N+q$.
Multiplication by $(t+(1-t))^{q-j}=1$ gives
\[
 t^{N+j}(1-t)^N
 =\sum_{\ell=0}^{q-j}\binom{q-j}{\ell}
 t^{N+j+\ell}(1-t)^{N+q-j-\ell}.
\]
Thus
\[
 t^N(1-t)^NQ(t)=\sum_{k=0}^n b_k t^k(1-t)^{n-k},
\]
where only $N\le k\le N+q$ can have $b_k\ne0$, and
\[
 |b_k|\le C_Q:=\sum_{j=0}^q |a_j|2^{q-j}.
\]
This bound does not depend on $N$. On the other hand, for $0\le j\le q$,
\[
 \binom{2N+q}{N+j}\ge\binom{2N}{N}
 =\prod_{h=1}^N\frac{N+h}{h}\ge2^N.
\]
The first inequality above can also be seen by adjoining a fixed set
of $j$ of the $q$ extra elements to each $N$-element subset of a
$2N$-element set. Choose $N$ so that $2^N\ge C_Q$. Put
$b_k^+=\max\{b_k,0\}$ and $b_k^-=\max\{-b_k,0\}$.
Both coefficient vectors satisfy the bounds in
\eqref{eq:Bernstein-word}.
Choose the leaf labels of $U,V$ to realize these two vectors. Their
polynomial difference is \eqref{eq:polynomial-realization}.
\end{proof}

\subsection{Positive differences of word values}

Fix $a<b$, let $A$ be the generated subalgebra, and restrict the fixed
affine embedding to $A$, using $a$ as base point. A \emph{positive word difference} means an element
$j(v)-j(u)$ with $u<v$ in $A$. This terminology refers to the real order
on $A$; it does not introduce an order on all of $M$.

\begin{lem}
\label{lem:no-positive-cancellation}
If $u_i<v_i$ in $A$ for $1\le i\le m$, then
\begin{equation}
 \sum_{i=1}^m(j(v_i)-j(u_i))\ne0.
 \label{eq:no-positive-cancellation}
\end{equation}
\end{lem}

\begin{proof}
The reason for choosing leaves with equal coefficients is already visible
at depth two. The middle two leaves of $(z_1*z_2)*(z_3*z_4)$ both have
coefficient map $\alpha\beta$. Thus the image difference between
$(a*v_1)*(v_2*a)$ and $(a*u_1)*(u_2*a)$ is
\[
 (\alpha\beta)\bigl((j(v_1)-j(u_1))+(j(v_2)-j(u_2))\bigr).
\]
The second evaluation is smaller than the first by strict increase.

For the general case choose $n,k$ with $\binom nk\ge m$. In a complete
binary expression of depth $n$, select $m$ leaves with coefficient
$\alpha^k\beta^{n-k}$. Fill the other leaves with $a$. Compare the
evaluation with $u_i$ at the selected leaves to the one with $v_i$ there.
Strict increase in every leaf makes the latter value larger, whereas their
image difference is
\[
 (\alpha^k\beta^{n-k})\left(\sum_{i=1}^m(j(v_i)-j(u_i))\right).
\]
If the sum were zero, the two images would agree, contradicting
injectivity of the embedding.
\end{proof}

\subsection{Positive coefficients in the binary-tree basis}

The basis $t^k(1-t)^{n-k}$ comes from the paths in a complete binary tree.
Its usefulness here is algebraic: after substituting $\alpha$ for $t$,
each basis term becomes a composition of the two coefficient maps.
We need an exact positive-coefficient expansion, rather than an
approximation of a function.
The following is the interval form of P\'olya's positivity theorem;
see, for example, \cite{PowersReznick2001}. We include a finite proof
adapted to the present basis.

\begin{lem}
\label{lem:positive-tree-expansion}
If $R\in\Z[t]$ is strictly positive on $[0,1]$, there are $n\ge0$ and
positive integers $C,c_0,\ldots,c_n$ such that
\begin{equation}
 C\cdot R(t)=\sum_{k=0}^n c_k t^k(1-t)^{n-k}.
 \label{eq:positive-Bernstein-expansion}
\end{equation}
\end{lem}

\begin{proof}
Write $R(t)=\sum_{j=0}^q a_jt^j$ and put
\[
 \delta:=\min_{t\in[0,1]}R(t)>0.
\]
We first write $R$ in the required basis. We then compare its coefficients,
after division by $\binom nk$, with the positive values $R(k/n)$.
An explicit bound for the difference will allow us to choose $n$ so that
all coefficients are positive.

\medskip\noindent
\emph{Step 1: an exact formula for the coefficients.}
Fix an integer $n\ge q$. For $0\le j\le q$, the binomial theorem gives
\[
 t^j=t^j\bigl(t+(1-t)\bigr)^{n-j}
 =\sum_{k=j}^n\binom{n-j}{k-j}t^k(1-t)^{n-k}.
\]
For $j\le k\le n$ we have
\[
 \binom{n-j}{k-j}
 =\binom nk\frac{\binom kj}{\binom nj}.
\]
Substituting these identities into $R$ and collecting the coefficients
therefore gives
\begin{equation}
 \begin{aligned}
 R(t)&=\sum_{k=0}^n b_{n,k}\binom nk t^k(1-t)^{n-k},\\
 b_{n,k}&=\sum_{j=0}^q a_j\frac{\binom kj}{\binom nj}
 \qquad(0\le k\le n),
 \end{aligned}
 \label{eq:tree-basis-coefficients}
\end{equation}
where $\binom kj=0$ when $k<j$.
In particular, every $b_{n,k}$ is rational.

For example, if $R(t)=a_0+a_1t+a_2t^2$ and $n\ge2$, this formula reads
\[
 b_{n,k}
 =a_0+a_1\frac{k}{n}
   +a_2\frac{k(k-1)}{n(n-1)}.
\]
The corresponding value $R(k/n)$ has $(k/n)^2$ in the last term.
The next step bounds this type of difference for every degree $j$.

\medskip\noindent
\emph{Step 2: compare the coefficient factors with powers of $k/n$.}
We claim that, whenever $n\ge\max\{1,2q\}$,
\begin{equation}
 \left|
 \frac{\binom kj}{\binom nj}-\left(\frac{k}{n}\right)^j
 \right|
 \le\frac{j(j-1)}{n}
 \qquad(0\le k\le n,\ 0\le j\le q).
 \label{eq:finite-factorial-bound}
\end{equation}
For $j=0,1$, the two expressions agree.
Suppose that $j\ge2$. If $k<j$, the binomial quotient is zero and
\[
 0\le\left(\frac{k}{n}\right)^j
 \le\frac{k}{n}\le\frac{j-1}{n}
 \le\frac{j(j-1)}{n}.
\]
It remains to consider $k\ge j$. In this case
\[
 \frac{\binom kj}{\binom nj}
 =\prod_{h=0}^{j-1}\frac{k-h}{n-h}.
\]
We compare each factor with $k/n$. Since $h<j\le q$ and $n\ge2q$,
\[
 0\le\frac{k}{n}-\frac{k-h}{n-h}
 =\frac{h(n-k)}{n(n-h)}
 \le\frac{2h}{n}.
\]
All factors lie in $[0,1]$. Replacing one factor in a product by another
number in $[0,1]$ changes the product by at most the change in that
factor, since the product of the remaining factors is at most $1$.
Replacing the $j$ factors one at a time thus yields
\[
 0\le
 \left(\frac{k}{n}\right)^j-\frac{\binom kj}{\binom nj}
 \le
 \sum_{h=0}^{j-1}
 \left(\frac{k}{n}-\frac{k-h}{n-h}\right)
 \le\sum_{h=0}^{j-1}\frac{2h}{n}
 =\frac{j(j-1)}{n}.
\]
This proves \eqref{eq:finite-factorial-bound} in all cases.

\medskip\noindent
\emph{Step 3: choose the degree and clear the denominators.}
Set
\[
 E_R:=\sum_{j=0}^q |a_j|j(j-1).
\]
By \eqref{eq:tree-basis-coefficients} and
\eqref{eq:finite-factorial-bound}, for every $0\le k\le n$,
\[
 \begin{aligned}
 |b_{n,k}-R(k/n)|
 &\le
 \sum_{j=0}^q |a_j|
 \left|
 \frac{\binom kj}{\binom nj}
 -\left(\frac{k}{n}\right)^j
 \right|\\
 &\le\frac{E_R}{n}.
 \end{aligned}
\]
The number $E_R$ depends only on $R$, so we may choose an integer
$n\ge\max\{1,2q\}$ with $n>E_R/\delta$.
For this choice, all coefficients satisfy
\[
 b_{n,k}\ge R(k/n)-\frac{E_R}{n}
 \ge\delta-\frac{E_R}{n}>0.
\]
Finally, choose a positive integer $C$ divisible by the denominators
of the finitely many rational numbers $b_{n,0},\ldots,b_{n,n}$.
Then
\[
 c_k:=C\binom nk b_{n,k}\qquad(0\le k\le n)
\]
are positive integers. Multiplying
\eqref{eq:tree-basis-coefficients} by $C$ gives
\eqref{eq:positive-Bernstein-expansion}, as required.
\end{proof}

\begin{cor}
\label{cor:positive-polynomial}
Let $R\in\Z[t]$ be strictly positive on $[0,1]$, and let $u<v$ in $A$.
Then $R(\alpha)(j(v)-j(u))\ne0$. More generally, this element cannot
cancel a finite sum of positive word differences.
\end{cor}

\begin{proof}
Apply the same translation to the two points $u<v$. The map $z\mapsto z*a$
changes their image difference to $\alpha(j(v)-j(u))$, and $z\mapsto a*z$
changes it to $\beta(j(v)-j(u))$. Both translations are strictly increasing.
Thus every
$(\alpha^k\beta^{n-k})(j(v)-j(u))$ is again a positive word difference.

By Lemma~\ref{lem:positive-tree-expansion}, a positive integer multiple of
$R(\alpha)(j(v)-j(u))$ is a sum of such differences. If this element,
alone or together with a finite sum of positive word differences, were
zero, clearing the same positive integer would contradict
Lemma~\ref{lem:no-positive-cancellation}.
\end{proof}

\subsection{The Euclidean algorithm with a positive multiplier}

For relatively prime polynomials $f,g\in\Q[t]$, there are polynomials
$p,q\in\Q[t]$ satisfying
\[
 p(t)f(t)+q(t)g(t)=1.
\]
This follows from the Euclidean algorithm\footnote{The displayed relation
is usually called the polynomial B\'ezout identity.}; see
\cite[1.8]{MilneFT}. The next lemma adjusts such an identity so that
the polynomial multiplying $H$ is positive on the whole interval $[0,1]$.
Positivity only at $r$ would not suffice to apply
Corollary~\ref{cor:positive-polynomial}.

\begin{lem}
\label{lem:positive-bezout}
Let $m\in\Z[t]$ be irreducible over $\Q$, with a unique zero
$r\in[0,1]$. If $H\in\Z[t]$ and $H(r)>0$, there are
$A_1,B_1\in\Z[t]$ and an integer $C>0$ such that
\begin{equation}
 A_1(t)>0\quad(0\le t\le1),\qquad A_1H+B_1m=C.
 \label{eq:positive-Bezout}
\end{equation}
\end{lem}

\begin{proof}
Since $m$ is irreducible and $H(r)>0$, the polynomial $m$ does not
divide $H$. Hence $H,m$ are relatively prime in $\Q[t]$.
The identity stated above gives $A_0,B_0\in\Q[t]$ with
\[
 A_0H+B_0m=1.
\]
At $t=r$ this identity gives
$A_0(r)=1/H(r)>0$.

If $A_0$ is positive on $[0,1]$, put $A=A_0$ and $B=B_0$.
Otherwise the compact set
\[
 E:=\{t\in[0,1]:A_0(t)\le0\}
\]
is nonempty and does not contain $r$. Since $r$ is the only zero of $m$
in $[0,1]$, the number $\eta:=\min_E m(t)^2$ is positive.
Choose an integer $K>\max_{[0,1]}|A_0|/\eta$ and put
\[
 A:=A_0+K m^2,\qquad B:=B_0-KmH.
\]
On $E$ the added term exceeds $|A_0|$, and outside $E$ the
polynomial $A_0$ is already positive. Thus $A>0$ on $[0,1]$.
The adjustment preserves the identity, since
\[
 AH+Bm=(A_0+Km^2)H+(B_0-KmH)m=1.
\]
In either case, choose a positive integer $C$ clearing the
denominators of $A,B$, and set $A_1=CA$, $B_1=CB$.
Then $A_1,B_1\in\Z[t]$, $A_1>0$ on $[0,1]$, and
$A_1H+B_1m=C$, as required.
\end{proof}

\section{Unique-root word identities}

The polynomial tools allow us to construct the weighted representation
from a single word coincidence.

Lemma~\ref{lem:word-coincidence-subalgebra}(ii) already extends the word
identity to the subalgebra generated by the original pair.
The global conclusion of Corollary~\ref{cor:polynomial-globalization}
additionally requires that the classes be intervals. The root hypothesis
in the next theorem supplies the local representation needed to establish
this property.

\begin{thm}
\label{thm:unique-root}
Let $a<b$ in $I$, and let $U,V$ be binary words such that
\begin{equation}
 U(a,b)=V(a,b).
 \label{eq:unique-root-equality}
\end{equation}
Suppose that $Q:=P_U-P_V\ne0$ has exactly one zero $r$ in $(0,1)$,
and that this zero is simple. Then there is a continuous strictly
increasing bijection $\varphi:[a,b]\to[0,1]$ such that
\begin{equation}
 \varphi(F(u,v))=r\varphi(u)+(1-r)\varphi(v)
 \qquad(u,v\in[a,b]).
 \label{eq:unique-root-representation}
\end{equation}
Moreover, $U(x,y)=V(x,y)$ for all $x,y\in I$, and there is a continuous
strictly increasing function $\psi:I\to\R$ with
\[
 F(u,v)=\psi^{-1}\bigl(r\psi(u)+(1-r)\psi(v)\bigr)
 \qquad(u,v\in I).
\]
\end{thm}

\begin{proof}
Use the fixed affine embedding of $(I,*)$, take $a$ as base point, and put
$d=j(b)$. Let $A$ be the subalgebra generated by $a,b$.
We retain the distinction between real points and their images from
Section~2; products of $\alpha,\beta$ below are compositions of maps.

\medskip\noindent
\emph{Step 1: reduce the given equality to the minimal polynomial.}
Let $m\in\Z[t]$ be the primitive minimal polynomial of $r$.
Gauss's lemma gives
\begin{equation}
 Q(t)=t^p(1-t)^q m(t)R(t),
 \label{eq:unique-root-factorization}
\end{equation}
where $p,q\ge0$ and $R\in\Z[t]$ has no zero on $[0,1]$.
Indeed, remove the endpoint factors and then $m$; an additional interior
zero would contradict the uniqueness of $r$, and an additional factor
vanishing at $r$ would contradict its simplicity. Choose
$R_+\in\{R,-R\}$ positive on $[0,1]$.

By Lemma~\ref{lem:word-evaluation}, the assumed equality is
$Q(\alpha)(d)=0$. Injectivity of $\alpha,\beta$ allows the endpoint
factors to be cancelled, and hence
\begin{equation}
 R_+(\alpha)\bigl(m(\alpha)(d)\bigr)=0.
 \label{eq:cofactor-relation}
\end{equation}
To use Corollary~\ref{cor:positive-polynomial}, we first realize a
multiple of $m(\alpha)(d)$ as a difference of word values.
Lemma~\ref{lem:polynomial-realization} gives words $W,Z$ and $N\ge0$ with
$P_W-P_Z=t^N(1-t)^Nm$. Therefore
\[
 e:=j(W(a,b))-j(Z(a,b))=(\alpha^N\beta^N)\bigl(m(\alpha)(d)\bigr).
\]
Equation \eqref{eq:cofactor-relation} gives $R_+(\alpha)(e)=0$.
If $e\ne0$, either $e$ or $-e$ is a positive word difference, contrary
to Corollary~\ref{cor:positive-polynomial}. Thus $e=0$; cancelling
$\alpha^N\beta^N$ gives
\begin{equation}
 m(\alpha)(d)=0.
 \label{eq:minimal-polynomial-relation}
\end{equation}

\medskip\noindent
\emph{Step 2: equality and order are detected by evaluation at $r$.}
Here $W,Z$ denote arbitrary binary words. If $P_W(r)=P_Z(r)$,
Gauss's lemma gives $m\mid P_W-P_Z$ in $\Z[t]$. Equation
\eqref{eq:minimal-polynomial-relation} and
Lemma~\ref{lem:word-evaluation} then imply
\begin{equation}
 P_W(r)=P_Z(r)\quad\Longrightarrow\quad W(a,b)=Z(a,b).
 \label{eq:unique-root-consistency}
\end{equation}

Suppose next that $H:=P_Z-P_W$ satisfies $H(r)>0$.
The polynomial $m$ has no other zero in $[0,1]$: any other interior zero
would be a zero of $Q$, and an irreducible polynomial with a root in
$(0,1)$ cannot vanish at $0$ or $1$.
Lemma~\ref{lem:positive-bezout} supplies
$A_1>0$ on $[0,1]$, $B_1\in\Z[t]$, and $C\in\N$, $C>0$, with
$A_1H+B_1m=C$. Put
$h=j(Z(a,b))-j(W(a,b))=H(\alpha)(d)$. Applying the identity to $d$ gives
\begin{equation}
 A_1(\alpha)(h)=Cd.
 \label{eq:ordered-Bezout}
\end{equation}
The right-hand side is a sum of $C$ copies of the positive word difference
$d=j(b)-j(a)$, so it is nonzero by
Lemma~\ref{lem:no-positive-cancellation}. Thus $h\ne0$.
If $Z(a,b)<W(a,b)$, then $-h$ is a positive word difference and
\[
 A_1(\alpha)(-h)+Cd=0,
\]
contradicting Corollary~\ref{cor:positive-polynomial}.
We conclude that
\begin{equation}
 P_W(r)<P_Z(r)\quad\Longrightarrow\quad W(a,b)<Z(a,b).
 \label{eq:unique-root-order}
\end{equation}

\medskip\noindent
\emph{Step 3: construct the increasing function and complete it.}
Define
\begin{equation}
 S_r:=\{P_W(r):W\in\mathcal W\},\qquad
 f_r(P_W(r)):=W(a,b).
 \label{eq:unique-root-coding}
\end{equation}
Equations \eqref{eq:unique-root-consistency} and
\eqref{eq:unique-root-order} make $f_r$ well-defined and strictly
increasing. There are countably many finite binary words, so $S_r$ is
countable. By the recursion for $P_W$, the set $S_r\subset[0,1]$ contains
$0,1$ and is closed under $(s,t)\mapsto rs+(1-r)t$.
Every word value belongs to $[a,b]$ by Lemma~\ref{lem:strict-mean}, and
$f_r(0)=a$, $f_r(1)=b$.

To check density of $S_r$, start from the endpoints $0,1$ and repeatedly
insert $rs+(1-r)t$ and $rt+(1-r)s$ between adjacent points $s<t$.
After deleting a possible repetition when $r=1/2$ and ordering the
inserted points, every new gap is at most
$\rho(t-s)$, where $\rho=\max\{r,1-r\}<1$.
After $n$ subdivisions all gaps are at most $\rho^n$. Given any open
subinterval of $[0,1]$, choose $n$ with $\rho^n$ smaller than its length;
the interval must then contain a subdivision point. Thus $S_r$ is dense.

The word recursion and well-definedness of $f_r$ give
\[
 F(f_r(s),f_r(t))=f_r\bigl(rs+(1-r)t\bigr)
 \qquad(s,t\in S_r).
\]
All hypotheses of Theorem~\ref{thm:dense-domain}
(see also \cite[Proposition~7.1]{KissShulman}) now hold. That theorem proves density
of the range and extends $f_r$ to a continuous strictly increasing
bijection $\widetilde f_r:[0,1]\to[a,b]$.
Set $\varphi=\widetilde f_r^{-1}$ to obtain
\eqref{eq:unique-root-representation}.

\medskip\noindent
\emph{Step 4: apply the globalization principle.}
The first three steps have established the local conclusion for any pair
at which the word coincidence holds. Define
\[
 x\sim y\quad\Longleftrightarrow\quad U(x,y)=V(x,y).
\]
By Lemma~\ref{lem:word-coincidence-subalgebra}(ii), this is the relation
$\sim_Q$. To verify the interval hypothesis of
Corollary~\ref{cor:polynomial-globalization}, let $x<y$ with $x\sim y$.
The local conclusion of Steps~1--3 gives the weighted representation on
$[x,y]$ with weight $r$. Evaluating the two words in this representation
and using $P_U(r)=P_V(r)$ yields
\[
 U(u,v)=V(u,v)\qquad(u,v\in[x,y]).
\]
Thus all points of $[x,y]$ belong to one equivalence class.
Each class is therefore an interval, and the original coincidence gives
$a\sim b$.
Corollary~\ref{cor:polynomial-globalization} now gives
$U(x,y)=V(x,y)$ for all $x,y\in I$.
The local conclusion gives continuity on the square of every compact
subinterval. Such squares contain relative neighborhoods of all points of
$I^2$, so $F$ is continuous on $I^2$. The continuous bisymmetric
representation theorem \cite[Section~2]{Aczel1948} gives a global weighted quasi-arithmetic representation
with some weight $s\in(0,1)$. Evaluating the global word identity in that
representation gives $Q(s)=0$; uniqueness of the interior zero forces
$s=r$.
\end{proof}

The theorem also shows exactly how a suitable algebraic number can be encoded
by one word identity.

\begin{cor}
\label{cor:algebraic-construction}
Let $r\in(0,1)$ be algebraic, and suppose that its minimal polynomial has no
other zero in $(0,1)$.  There are binary words $U_r,V_r$ such that, whenever
$a<b$ and
\[
 U_r(a,b)=V_r(a,b),
\]
the operation $F$ has a global weighted quasi-arithmetic representation
with weight $r$, whose normalized restriction to $[a,b]$ is
\eqref{eq:unique-root-representation}.
\end{cor}

\begin{proof}
Let $m_r\in\Z[t]$ be the primitive minimal polynomial of $r$.  By
Lemma~\ref{lem:polynomial-realization}, there are words $U_r,V_r$ and an
$N\ge0$ such that
\[
 P_{U_r}(t)-P_{V_r}(t)=t^N(1-t)^Nm_r(t).
\]
Minimal polynomials over $\Q$ have only simple zeros.  The displayed
difference therefore has $r$ as its unique simple zero in $(0,1)$, and
Theorem~\ref{thm:unique-root} applies.
\end{proof}

\begin{proof}[Proof of Theorem~\ref{thm:intro-main}]
The forward implication, including the global conclusion, follows from Theorem~\ref{thm:unique-root}.

For the converse, let $m$ denote the primitive minimal polynomial of $r$. The construction in the proof of Corollary~\ref{cor:algebraic-construction}  gives two binary words whose polynomial difference is $m$ multiplied by powers of $t$ and of $(1-t)$.

These factors are nonzero between 0 and 1, so the difference has exactly the same interior zeros, with the same multiplicities, as $m$. Minimal polynomials over the rationals have only simple zeros, and by assumption $r$ is the only zero of $m$ between 0 and 1. Hence the constructed difference has all the required properties.
\end{proof}

\subsection{Examples: depth two and beyond}

Reflexivity allows every word of depth at most two to be written as a
complete depth-two word. By \eqref{eq:Bernstein-word}, every difference
of two such word polynomials has the form
\[
 Q(t)=c_0(1-t)^2+c_1t(1-t)+c_2t^2,
 \quad c_0,c_2\in\{-1,0,1\},\quad c_1\in\{-2,-1,0,1,2\}.
\]
If $Q\equiv0$, Lemma~\ref{lem:word-evaluation} shows that the equality
holds identically and gives no additional restriction on $F$.
If $Q$ has no zero in $(0,1)$, a coincidence at distinct points is
impossible: in the affine relation, cancel the powers of the injective
maps $\alpha,\beta$ corresponding to endpoint factors. The remaining
polynomial has a fixed strict sign on $[0,1]$, contrary to
Corollary~\ref{cor:positive-polynomial}.

The remaining possibilities can be found without listing all pairs of
words, by considering $Q(0)=c_0$ and $Q(1)=c_2$.
\begin{itemize}
\item If $c_0c_2<0$, interchange the two sides so that $c_0=-1$ and
$c_2=1$. Then $Q(t)=-c_1t^2+(2+c_1)t-1$.
The choices $c_1=-2,-1,0$ give $2t^2-1$, $t^2+t-1$, and $2t-1$;
the choices $c_1=1,2$ give the reflected roots $r\mapsto1-r$.
\item If $c_0c_2=0$, removing endpoint factors leaves a polynomial of
degree at most one. Its only possible interior zeros are $1/3$, $1/2$,
and $2/3$.
\item If $c_0c_2>0$, change the overall sign so that $c_0=c_2=1$.
Now $Q(t)=1+(c_1-2)t(1-t)$, whose only possible interior zero occurs
for $c_1=-2$, when $Q(t)=(2t-1)^2$.
\end{itemize}

The last, repeated-root case is, up to interchanging the sides,
\[
 (y*x)*(x*y)=(x*y)*(y*x).
\]
Suppose it holds at $a<b$, and put $p=a*b$, $q=b*a$.
If $p=q$, then $a*b=b*a$; otherwise the displayed equality gives
$p*q=q*p$ for distinct $p,q$. In either case,
\cite[Theorem~3]{BuraiKissSzokol2023} gives a global quasi-arithmetic
representation with weight $1/2$.

Every other remaining difference has exactly one simple zero in $(0,1)$,
so Theorem~\ref{thm:unique-root} applies. Apart from the symmetric weight
$1/2$, there are three cases up to replacing $*$ by the opposite
operation $u\circ v:=v*u$, which replaces $r$ by $1-r$.
Representatives are
\begin{align}
 (y*x)*x&=(x*y)*(y*x),
 &Q(t)&=t(3t-2), &r&=\frac23,                  \label{eq:example-23}\\
 (y*x)*x&=x*y,
 &Q(t)&=t^2+t-1, &r&=\frac{\sqrt5-1}{2},       \label{eq:example-golden}\\
 (y*x)*x&=(x*y)*y,
 &Q(t)&=2t^2-1, &r&=\frac1{\sqrt2}.            \label{eq:example-sqrt}
\end{align}
The opposite operation supplies the other three nonsymmetric weights
$1/3$, $(3-\sqrt5)/2$, and $1-1/\sqrt2$.

For $r=2/3$, the set of polynomial values $S_r$ is the set $D_3$ of
triadic rationals by Lemma~\ref{lem:canonical-coordinate}.
Reversing the arguments at every node of a word gives $S_{1-r}=S_r$,
and hence the same conclusion for $r=1/3$.
For the other nonsymmetric weights, $D_3$ is replaced by the dense set
obtained from $0,1$ by repeatedly applying
$(s,t)\mapsto rs+(1-r)t$; its elements lie in $\Q(r)\cap[0,1]$.
No separate recursive order calculation is needed, because Step~2 of the
proof of Theorem~\ref{thm:unique-root} establishes the order of all word
values at once.

These depth-two cases are only the first instances of the method.
A direct generalization of the $2/3$ identity is obtained by defining
formal words
\[
 U_1:=Y*X,\qquad V_1:=X*Y,\qquad
 U_{n+1}:=U_n*X,\qquad V_{n+1}:=V_n*U_n\quad(n\ge1).
\]
Induction using the word-polynomial recursion gives
\[
 P_{U_n}(t)=t^n,\qquad
 P_{V_n}(t)=n t^{n-1}(1-t).
\]
Consequently
\begin{equation}
 P_{U_n}(t)-P_{V_n}(t)
 =t^{n-1}\bigl((n+1)t-n\bigr).
 \label{eq:rational-family}
\end{equation}
Its only zero in $(0,1)$ is $r_n=n/(n+1)$, and this zero is simple.
Thus any coincidence $U_n(a,b)=V_n(a,b)$ with $a<b$ forces a global
weighted quasi-arithmetic representation with weight $r_n$, by
Theorem~\ref{thm:unique-root}.
For $n=2$ this is exactly \eqref{eq:model-identity-intro}; the cases
$n\ge3$ give explicit higher-depth identities, beginning with weight
$3/4$.

\section{Further directions}

Theorem~\ref{thm:unique-root} treats the case in which one word equality
selects a unique simple interior root. After the endpoint factors and its
minimal polynomial are removed, the remaining factor has a fixed sign on
$[0,1]$. Two algebraic cases lie beyond this argument. If the
word-polynomial difference has several distinct zeros in $(0,1)$, the
equality alone does not choose among them. If its only interior zero is
repeated, the affine calculation initially gives a relation of the form
\[
 m_r(\alpha)^k(d)=0\qquad(k>1).
\]
The proof above does not directly reduce this to $m_r(\alpha)(d)=0$.
That reduction requires a further argument. When there are several distinct
interior zeros, additional word information may help distinguish the possible
weights.

Polynomial equalities are necessarily insensitive to a transcendental
weight.  Indeed, for transcendental $r$,
\[
 P_U(r)=P_V(r)\quad\Longrightarrow\quad P_U=P_V.
\]
The order of word values, rather than their nontrivial equality, must then
carry the parameter.  In an affine model it determines the signs of all
differences $P_U(r)-P_V(r)$.  Lemma~\ref{lem:polynomial-realization} shows
that the complete sign pattern distinguishes any two points of $(0,1)$, so
it is a genuine analogue of a Dedekind cut.  Turning this sign information
intrinsically into a consistent ordered coding is the natural continuation
of the constructive method.  Once such a coding is obtained, the same set
$S_r$ and the same dense-domain theorem give the full weighted
quasi-arithmetic representation.

\section*{Acknowledgements.}
The project was supported by the Hungarian National Research, Development
and Innovation Office grants K146922, STARTING 150576.
P. Szokol was supported by the HUN-REN Hungarian Research Network as a member of HUN-REN-UD Equations, Functions, Curves and their
Applications Research Group.


\section*{Declaration of generative AI and AI-assisted technologies}
During the preparation of this work, ChatGPT 6 Sol (OpenAI)
was used to assist with editing the manuscript. In particular, we used to search
for the positive-coefficient representation needed in
Lemma~\ref{lem:positive-tree-expansion}.
It helped identify the relevant result as the interval form of
P\'olya's positivity theorem and assisted in developing and
presenting the proof of that lemma.
The authors reviewed and edited the material as needed and take
full responsibility for the content and integrity of the publication.


\section*{Appendix}

\subsection*{A. The binary dense-domain no-gap argument}
\setcounter{equation}{0}
\renewcommand{\theequation}{A.\arabic{equation}}
\renewcommand{\theHequation}{A.\arabic{equation}}

The function built from word values is initially known only on a dense
set of real parameters. Density of that parameter set does not by itself
imply density of its increasing image: the function could have a jump.
The exact averaging identity rules out this possibility. A hypothetical
gap would give a nonempty interval for each of uncountably many parameters,
and strict increase would force these intervals to be disjoint. This is
the countability argument needed to complete the interval.

We prove Theorem~\ref{thm:dense-domain} in the precise form
used in this paper.  It is the nonsymmetric binary version of the no-gap
argument used in \cite{BuraiKissSzokol2021}, and also the binary specialization
of \cite[Proposition~7.1]{KissShulman}.

Let $J=[a,b]$, let $0<r<1$, and assume that $S$, $H$, and $f$ satisfy the
hypotheses of Theorem~\ref{thm:dense-domain}.  Define
\begin{equation}
 \widehat f(t):=\sup\{f(s):s\in S,\ s\le t\}
 \qquad(0\le t\le1).
 \label{eq:monotone-extension}
\end{equation}
Then $\widehat f$ extends $f$ and is increasing.  In fact it is strictly
increasing: if $t_1<t_2$, choose $s_1,s_2\in S$ with
$t_1<s_1<s_2<t_2$.  Then
\[
 \widehat f(t_1)\le f(s_1)<f(s_2)\le\widehat f(t_2).
\]
Moreover $\widehat f([0,1])$ is contained in the closure of $f(S)$.
Consequently that closure is an uncountable compact subset of $J$.  A compact
subset of the real line has at most countably many isolated points and at
most countably many points that are accumulation points from only one side.
The latter are endpoints of complementary open intervals; there are at most
countably many such intervals because each contains a rational number.
It follows that the closure of $f(S)$ has uncountably many points that are
approached by points of $f(S)$ from both sides.

We claim that $f(S)$ is dense in $J$. Suppose otherwise. Choose
$a<c_0<d_0<b$ with $(c_0,d_0)\cap f(S)=\varnothing$, and put
\[
 z:=\sup\{s\in S:f(s)\le c_0\}.
\]
The set in this supremum contains $0$ and does not contain $1$.
For $s\in S$ with $s<z$ we have $f(s)\le c_0$; for $s>z$ we have
$f(s)\ge d_0$. Consequently the numbers
\begin{equation}
 X=\sup\bigl(\{a\}\cup\{f(s):s\in S,\ s<z\}\bigr),
 \quad
 Y=\inf\bigl(\{b\}\cup\{f(s):s\in S,\ s>z\}\bigr)
 \label{eq:assumed-gap}
\end{equation}
satisfy $X\le c_0<d_0\le Y$. The endpoint choices $z=0,1$ include a
possible jump next to $a$ or $b$. If $z\in S$, the image $f(z)$ is excluded
from the two sets defining $X,Y$ and does not affect this inequality.
Choose two distinct two-sided accumulation points $\xi<\eta$ of $f(S)$, and
then choose $p,q\in S$ so that
\[
 \xi<f(p)<f(q)<\eta.
\]
Because $S$ is dense and contains $0,1$, we may choose $d,D\in S$ with
$d\le z\le D$ and $d<D$, as close to $z$ as necessary, so that
\[
 f(d)\le X<Y\le f(D)
\]
(take $d=0$ if $z=0$, and $D=1$ if $z=1$)
and
\begin{equation}
 D-d<\frac{1-r}{r}(q-p).
 \label{eq:gap-parameter-choice}
\end{equation}
The two numbers
\[
 A_1:=rD+(1-r)p,
 \qquad
 A_2:=rd+(1-r)q
\]
belong to $S$, and \eqref{eq:gap-parameter-choice} says precisely that
$A_1<A_2$.  Hence the strict increase of $f$ and the dense-domain identity
give
\[
 H(f(D),f(p))=f(A_1)<f(A_2)=H(f(d),f(q)).
\]
Together with strict increase of $H$ in each variable, this yields
\[
 H(f(d),f(p))
 <H(f(D),f(p))
 <H(f(d),f(q))
 <H(f(D),f(q)).
\]
Using $f(d)\le X<Y\le f(D)$ in the middle inequality, we obtain
\[
 H(Y,f(p))
 \le H(f(D),f(p))<H(f(d),f(q))
 \le H(X,f(q)),
\]
so $H(Y,f(p))<H(X,f(q))$.  This remains strict even in the endpoint
cases.  We then obtain
from $\xi<f(p)<f(q)<\eta$,
\begin{equation}
 H(Y,\xi)<H(X,\eta).
 \label{eq:separated-gaps}
\end{equation}

The same reasoning applies to every pair $\xi<\eta$ of two-sided
accumulation points of $f(S)$. Therefore the nonempty open intervals
\[
 \Phi_t:=\bigl(H(X,t),H(Y,t)\bigr)
\]
are pairwise disjoint as $t$ ranges over the uncountable set of such points:
if $t_1<t_2$, then \eqref{eq:separated-gaps} gives
$H(Y,t_1)<H(X,t_2)$. This is impossible, since every nonempty open interval
contains a rational number and the rational numbers are countable.  Thus the
assumed gap cannot exist, and $f(S)$ is dense in $J$.

It now follows from \eqref{eq:monotone-extension} that $\widehat f$ has no
jump: any jump would leave a nonempty open interval disjoint from $f(S)$.
Hence $\widehat f$ is a continuous strictly increasing bijection from
$[0,1]$ onto $J$.  Denote it by $\widetilde f$.

It remains only to extend the identity.  Given $u,v\in J$, put
$s=\widetilde f^{-1}(u)$ and $t=\widetilde f^{-1}(v)$.  Choose sequences
$s_n^-,s_n^+,t_n^-,t_n^+\in S$ such that
\[
 s_n^-\le s\le s_n^+,
 \qquad
 t_n^-\le t\le t_n^+,
\]
and all four sequences converge to the corresponding limit.  Monotonicity
and the identity on $S$ give
\begin{align*}
 \widetilde f\bigl(rs_n^-+(1-r)t_n^-\bigr)
 &\le H(u,v)\\
 &\le
 \widetilde f\bigl(rs_n^++(1-r)t_n^+\bigr).
\end{align*}
Letting $n\to\infty$ and using the continuity of $\widetilde f$, we conclude
that
\[
 H(u,v)=\widetilde f\bigl(rs+(1-r)t\bigr).
\]
The extension is unique because two continuous extensions agreeing on the
dense set $S$ agree everywhere. This proves all assertions of
Theorem~\ref{thm:dense-domain}.

\subsection*{B. Proof of the globalization principle}
\setcounter{equation}{0}
\renewcommand{\theequation}{B.\arabic{equation}}
\renewcommand{\theHequation}{B.\arabic{equation}}

We prove Proposition~\ref{prop:global-propagation}, following the
equivalence-class strategy of
\cite[Section~3.2, proof of Theorem~3]{BuraiKissSzokol2023}.
There are two separate points to establish. First, the classes are closed intervals, even though
closedness was not assumed. Second, one non-degenerate class limits the
number of classes to a countable set, which cannot partition a compact
interval into disjoint closed intervals.

\medskip\noindent
\emph{Step 1: order and cancellation on classes.}
By \eqref{eq:quotient-cancellation}, replacing either input by an
equivalent point does not change the class of the output. Hence $*$
defines an operation on classes. Since the classes are disjoint intervals,
they inherit the order of $K$:
$C<D$ means that every point of $C$ is less than every point of $D$.
Fix $x\in C$, $y\in D$, and $z$ in an arbitrary class. Strict increase gives
$x*z<y*z$, and \eqref{eq:quotient-cancellation} says their classes are
distinct. Thus the operation on classes is strictly increasing in each
variable. Reflexivity gives $C*C=C$, so
\[
 C<C*D<D\qquad(C<D).
\]

\medskip\noindent
\emph{Step 2: the classes are closed.}
Let $C$ be a class and suppose $c:=\sup C$ does not belong to $C$.
Compactness of $K$ ensures $c\in K$. Let $D$ be the class containing $c$;
then $C<D$. For $x\in C$, the point $x*c$ belongs to $C*D$, a class
strictly above $C$. Therefore $x*c\ge\sup C=c$. On the other hand,
the strict mean property gives $x*c<c$. This contradiction shows $c\in C$.
The argument at $\inf C$ is the same. Each class is consequently a closed
interval, possibly a singleton.

\medskip\noindent
\emph{Step 3: there are at most countably many classes.}
Let $C_0$ contain $a<b$. For every class $D$, choose $y\in D$.
The two points $a*y<b*y$ belong to $C_0*D$, so this class contains the
nonempty open interval $(a*y,b*y)$. Different $D$ give different classes
$C_0*D$ by cancellation on classes. These open intervals are disjoint.
Choose a rational number in each of them; their disjointness makes the
chosen rationals distinct. Since $\Q$ is countable, there can be only
countably many classes $D$.

\medskip\noindent
\emph{Step 4: a countable partition by closed intervals is impossible.}
Suppose there is more than one class, and remove from $K$ the relative
interior of every class. Denote the remaining set by $E$.
It is compact. It is nonempty, since otherwise the distinct, disjoint
relative interiors would separate the connected interval $K$ into
nonempty open sets. Each closed interval class contributes at most two
points to $E$, so $E$ is countable by Step~3.

The set $E$ has no isolated point. Indeed, if an interior point $z\in E$
were isolated, sufficiently small intervals immediately to its left and
right would consist of points in relative interiors of classes. Each
of these intervals is connected, so all its points lie in a single class.
Both neighboring classes are closed and hence contain $z$. They are
therefore the same class, making $z$ an interior point of that class,
contrary to $z\in E$. At an endpoint of $K$, the same one-sided argument
shows that the neighboring class contains a relative neighborhood of $z$.

Finally, a nonempty compact subset of $\R$ without isolated points is
uncountable. For completeness, choose two disjoint closed subintervals,
each with a point of $E$ in its interior, and within each choose two
smaller subintervals with the same property.
Continue with lengths at most $2^{-n}$ at stage $n$. This is possible
because no point of $E$ is isolated. Every infinite binary sequence
selects nested closed intervals whose intersection contains one point
of $E$ by compactness, and different sequences give different points.
This contradicts countability of $E$. Thus there is only one class.

\end{document}